\documentclass{elsart}
\usepackage{amssymb, amsmath, url,graphicx}
\journal{Journal of Algebra}

\def\calm{{\cal M}}
\def\3{\subset }
\def\4{\subseteq }
\def\<{\left<}
\def\>{\right>}
\def\vsp{\vspace*{1,5mm}\\ }

\def\bit{\begin{itemize}}
\def\eit{\end{itemize}}
\def\3{\subset }
\def\4{\subseteq }
\def\ov{\overline}

\def\0{\leqno}

\def\barr{\begin{array}}
\def\earr{\end{array}}
\def\dd{\displaystyle}

\def\Z{{\rlap{$\kern2pt{\rm Z}$}{\rm Z}\,}}
\def\bld#1#2{{\buildrel{#1}\over{#2}}}
\def\st#1#2{{\mathrel{\mathop{#2}\limits_{#1}}{}\!}}
\def\stb#1#2#3{{\st{{#1}}{\bld{{#2}}{#3}}{}\!}}
\def\xmare#1#2{\stb{#1}{#2}{\mbox{\Huge$\times$}}}

\def\frax{\dd\frac}

\begin{document}
\begin{frontmatter}

\title{{\rm Addendum}\\ Addendum to "Subgroup commutativity degrees of finite groups"\\ $\left[\text{J. Algebra 321 (2009), 2508-2520}\right]$}
\author{Marius T\u arn\u auceanu}
\ead{tarnauc@uaic.ro}
\address{Faculty of Mathematics,
\textquotedblleft Al.I. Cuza\textquotedblright\ University, Ia\c si, Romania}

\begin{abstract}
We indicate a natural generalization of the concept of subgroup
commutativity degree of a finite group and a list of open problems
on these new concepts.
\end{abstract}

\end{frontmatter}

\section{Introduction}

The starting point for our discussion is given by \cite{4}, where
the {\it subgroup commutativity degree} of a finite group $G$ has
been introduced and studied. This new quantity is defined by
$$\barr{lcl}
sd(G)&=&\frax1{|L(G)|^2}\,\left|\{(H,K)\in L(G)^2\mid
HK=KH\}\right|=\vsp &=&\frax1{|L(G)|^2}\,\left|\{(H,K)\in
L(G)^2\mid HK\in L(G)\}\right|\earr$$and measures the probability
that two subgroups of $G$ commute, or equivalently the probability
that the product of two subgroups of $G$ be a subgroup in $G$. It
was inspired by the well-known commutativity degree $d(G)$ of $G$.
Since for $d(G)$ there is a natural generalization, namely the
{\it relative commutativity degree} of $G$ (see \cite{2}), a
similar one can be introduced for $sd(G)$. So, we define the {\it
relative subgroup commutativity degree} of a subgroup $H$ of $G$
$$sd(H,G)=\frax1{|L(H)||L(G)|}\,\left|\{(H_1,G_1)\in L(H)\times
L(G)\mid H_1G_1=G_1H_1\}\right|,$$and, more generally, the {\it
relative subgroup commutativity degree} of two subgroups $H$ and
$K$ of $G$
$$sd(H,K)=\frax1{|L(H)||L(K)|}\,\left|\{(H_1,K_1)\in L(H)\times
L(K)\mid H_1K_1=K_1H_1\}\right|.$$It is obvious that
$sd(G)=sd(G,G)$, for any finite group $G$, and that the above two
notions also have a probabilistic significance. In the following
we shall focus on some basic properties of the relative subgroup
commutativity degree and on its connections with the classical
subgroup commutativity degree.

On the other hand, in the final section of \cite{4} some further
research directions and three open problems on subgroup
commutativity degrees have been indicated. Since this concept, as
well as its above generalizations are very new, we think that a
more large list of open problems can be useful.

\section{Relative subgroup commutativity degrees of finite groups}

Let $G$ be a finite group and $H$ be a subgroup of $G$. Then
$$0<sd(H,G)\le1.$$Obviously, the equality $sd(H,G)=1$ holds if and only if all subgroups of $H$ are permutable in
$G$, or equivalently if and only if $H$ is modular and subnormal
in $G$ (see Theorem 5.1.1 of \cite{3}).

If $H\neq G$, then the set $\{(H_1,G_1)\in L(H)\times L(G)\mid
H_1G_1=G_1H_1\}$ contains the union of the disjoint sets
$\{(H_1,G_1)\in L(H)^2\mid H_1G_1=G_1H_1\}$ and $\{(H_1,G)\mid
H_1\in L(H)\}$. This shows that $sd(H,G)$ and $sd(H)$ satisfy the
inequality
$$sd(H,G)\geq\frax{|L(H)|}{|L(G)|}\hspace{0,5mm}sd(H)+\frax1{|L(G)|}\hspace{0,5mm}.$$
In the following, for every $H_1\in L(H)$, we shall denote by
$C(H_1)$ the set of subgroups of $G$ which commute with $H_1$ and
by $I(H_1)$ the set of subgroups of $G$ strictly containing $H_1$.
One obtains
$$sd(H,G)=\frax1{|L(H)||L(G)|}\dd\sum_{H_1\in L(H)}|C(H_1)|.$$Clearly, $N(G)$ is contained in each set $C(H_1)$, which implies that
$$sd(H,G)\geq\frax{|N(G)|}{|L(G)|}\hspace{0,5mm}.$$Since $L(H_1)\cup I(H_1) \subseteq C(H_1)$, for all $H_1\in L(H)$, we also infer that
$$sd(H,G)\geq\frax1{|L(H)||L(G)|}\hspace{0,5mm}\left(\dd\sum_{H_1\in L(H)}|L(H_1)|+\dd\sum_{H_1\in
L(H)}|I(H_1)|\right).$$Moreover, if $H_1\in N(G)$, then we find
the following inequality between the relative subgroup
commutativity degrees of $H$ and of $H/H_1$:
$$sd(H,G)\geq\frax{|L(H/H_1)||L(G/H_1)|}{|L(H)||L(G)|}\hspace{1mm}sd(H/H_1,G/H_1).$$
We remark that the permutability of the subgroups $(H_1,G_1)\in
L(H)\times L(G)$ is equivalent to the permutability of the
subgroups $(H_1^x,G_1)\in L(H^x)\times L(G)$, for every $x\in G$.
This leads to the following proposition.

\noindent{\bf Proposition 2.1.} {\it Any two conjugate subgroups
of a finite group have the same relative subgroup commutativity
degree.}

In the following let $(G_i)_{i=\overline{1,k}}$ be a family of
finite groups having coprime orders. Then the subgroup lattice of
the direct product $\prod_{i=1}^k G_i$ is decomposable, that is
every subgroup $H$ of $\prod_{i=1}^k G_i$ can (uniquely) be
written as $H=\prod_{i=1}^k H_i$ with $H_i\leq G_i$, for all
$i=\overline{1,k}$. A result similar with Proposition 2.2 of
\cite{4} is now obtained for the relative subgroup commutativity
degree.

\noindent{\bf Proposition 2.2.} {\it Under the above hypotheses,
the following equality holds
$$sd(H,\prod_{i=1}^k G_i)=\prod_{i=1}^k sd(H_i,G_i).$$}Obviously,
the above formula can successfully be applied in the case of
finite nilpotent groups.

\noindent{\bf Corollary 2.3.} {\it Let $G$ be a finite nilpotent
group and $(G_i)_{i=\ov{1,k}}$ be the Sylow subgroups of $G$.
Then, for every subgroup $H$ of $G$, we have
$$sd(H,G)=\prod_{i=1}^k sd(H_i,G_i),$$where $H_i$, $i=1,2,...,k$, are the Sylow subgroups of $H$.
In particular, we infer that the computation of the relative
subgroup commutativity degrees of subgroups of finite nilpotent
groups is reduced to $p$-groups.}

Our next goal is to establish some connections between $sd(G)$ and
the relative subgroup commutativity degrees of the maximal
subgroups of $G$, say $M_0$, $M_1$, ..., $M_r$. Let $H\in L(G)$.
Then $L(G)=\{G\}\cup \big(\cup_{i=0}^rL(M_i)\big)$ and so
$C(H)=\{G\}\cup\big(\cup_{i=0}^r\calm_i(H)\big)$, where
$\calm_i(H)=\{K\in L(M_i) \mid HK=KH\}$, for any $i=\ov{0,r}$. By
applying the well-known Inclusion-Exclusion Principle, it follows
that
$$|C(H)|=1+\dd\sum_{s=0}^r(-1)^s\hspace{-5mm}\dd\sum_{0\leq i_0<i_1<...<i_s\leq
r}|\cap_{j=0}^s\calm_{i_j}(H)|\hspace{0,5mm}.$$Since
$$sd(G)=\frax1{|L(G)|^2}\dd\sum_{H\in L(G)}|C(H)|$$and
$$\dd\sum_{H\in L(G)}|\cap_{j=0}^s\calm_{i_j}(H)|=\dd\sum_{H\in L(G)}\hspace{-2mm}|\{K\in L(\cap_{j=0}^sM_{i_j})\mid
HK=KH\}|\hspace{-1mm}=$$
$$\hspace{30mm}=|L(G)|\hspace{0,5mm}|L(\cap_{j=0}^sM_{i_j})|\hspace{0,5mm}sd(\cap_{j=0}^sM_{i_j},G),$$we
have proved the following result.

\noindent{\bf Theorem 2.4.} {\it Let $G$ be a finite group and
$M_0, M_1, ..., M_r$ be the maximal subgroups of $G$. Then
$$sd(G)=\frax1{|L(G)|}\left(1+\dd\sum_{s=0}^r(-1)^s\hspace{-5mm}\dd\sum_{0\leq i_0<i_1<...<i_s\leq
r}\hspace{-3mm}|L(\cap_{j=0}^sM_{i_j})|\hspace{1mm}sd(\cap_{j=0}^sM_{i_j},G)\right).\0(1)$$}Clearly,
the above equality allows us to compute the subgroup
\,commutativity degree for all finite groups $G$ whose maximal
subgroup structure is known. We also remark that certain
supplementary assumptions on the maximal subgroups of $G$ can
simplify the right side of (1). One of them consists in asking
that the relative subgroup commutativity degree of any
intersection of at least two (distinct) maximal subgroups of $G$
be equal to 1. In this case $sd(G)$ will depend only on
$sd(M_i,G)$, $i=0,1,...,r$.

\noindent{\bf Corollary 2.5.} {\it Let $G$ be a finite group and
$M_0, M_1, ..., M_r$ be the maximal subgroups of $G$. If
$sd(\cap_{j=0}^sM_{i_j},G)=1$, for any $s=\ov{1,r}$ and $0\leq
i_0<i_1<...<i_s\leq r$, then we have
$$sd(G)=1-\frax1{|L(G)|}\dd\sum_{i=0}^r|L(M_i)|\hspace{0,5mm}(1-sd(M_i,G)),\0(2)$$or equivalently
$$sd(G)=1-\frax1{|L(G)|^2}\dd\sum_{i,j=0}^r|L(M_i)||L(M_j)|\hspace{0,5mm}(1-sd(M_i,M_j)).\0(3)$$}In
\cite{4}, the explicit value $sd(A_4)=16/25$ has been directly
computed. Since $A_4$ satisfies the supplementary condition in the
hypotheses of Corollary 2.5, this value can be also obtained by
using (2) or (3). The same thing cannot be said in the case of
$S_4$, for which we must apply the general formula (1).

\noindent{\bf Example 2.6.} It is well-known that $S_4$ possesses
eight maximal subgroups: $M_0=A_4$, $M_i\cong S_3$, for $1\leq i
\leq 4$, and $M_i\cong D_8$, for $5\leq i\leq 7$. By inspecting
$L(S_4)$, we infer that the intersections of any $s\geq5$ distinct
ma\-xi\-mal subgroups is trivial, while the intersections of
$s\leq4$ distinct maximal subgroups are isomorphic with
$\mathbb{Z}_2$, $\mathbb{Z}_3$, $\mathbb{Z}_2\times\mathbb{Z}_2$,
$S_3$, $D_8$ or $A_4$. Then (1) leads to
$sd(S_4)=\dd\frac{1}{30}\hspace{1mm}\big(13-24sd(\mathbb{Z}_2,S_4)-8sd(\mathbb{Z}_3,S_4)-18sd(\mathbb{Z}_2\times\mathbb{Z}_2,S_4)+
24sd(S_3,S_4)+30sd(D_8,S_4)+10sd(A_4,S_4)\big).$ We easily find:
$sd(\mathbb{Z}_2,S_4)=2/3$, $sd(\mathbb{Z}_3,S_4)=7/12$,
$sd({\mathbb{Z}_2{\times}\mathbb{Z}_2},S_4)=44/75$,
$sd(S_3,S_4)=4/9$, $sd(D_8,S_4)=37/75$ and $sd(A_4,S_4)=151/300$.
Hence $sd(S_4)=1841/4500$.

\section{Open problems}

\noindent{\bf Problem 3.1.} Let $G$ be a finite group and $H\in
L(G)$. Which are the connections between $sd(G)$ and the classical
commutativity degree $d(G)$, respectively between $sd(H,G)$ and
the classical relative commutativity degree $d(H,G)$?

\noindent{\bf Problem 3.2.} The relative subgroup commutativity
degrees can be obviously computed for finite groups whose subgroup
structure is precisely determined. An interesting example of such
groups is constituted by the finite groups with all Sylow
subgroups cyclic, the so-called ZM-groups. Such a group is of type
$${\rm ZM}(m,n,r)=< a, b \mid a^m = b^n = 1, \hspace{1mm}b^{-1} a b =
a^r>,$$where the triple $(m,n,r)$ satisfies the conditions ${\rm
gcd}(m,n) = {\rm gcd}(m, r-1) = 1$ and $r^n \equiv 1
\hspace{1mm}({\rm mod}\hspace{1mm}m)$. The subgroups of ${\rm
ZM}(m,n,r)$ have been completely described in \cite{1}. Set
$$L=\left\{(m_1,n_1,s)\in\mathbb{N}^3 \hspace{1mm}\mid\hspace{1mm}
m_1\mid m,\hspace{1mm} n_1\mid n,\hspace{1mm} s<m_1,\hspace{1mm}
m_1\mid s\frac{r^n-1}{r^{n_1}-1}\right\}.$$Then there is a
bijection between $L$ and $L({\rm ZM}(m,n,r))$, namely the
function that maps a triple $(m_1,n_1,s)\in L$ into the subgroup
$$H_{(m_1,n_1,s)}=\cup_{k=1}^{n/n_1}(b^{n_1}a^s)^k<a^{m_1}> \mbox{
of } \hspace{1mm}{\rm ZM}(m,n,r).$$Give an explicit formula for
$sd(H_{(m_1,n_1,s)},{\rm ZM}(m,n,r))$.

\noindent{\bf Problem 3.3.} It is clear that $sd(A_3,S_3)=1$. We
also have seen in Section 2 that $sd(A_4,S_4)=151/300$. These lead
to the following two natural asks: compute $sd(A_n,S_n)$, for an
arbitrary $n\geq5$, and the limit
$\dd\lim_{n\to\infty}sd(A_n,S_n)$.

\noindent{\bf Problem 3.4.} By using (1), for a finite group $G$
we are able to calculate $sd(G)$ whenever the structure of maximal
subgroups of $G$ and their relative subgroup commutativity degrees
are known. Is this true for other remarkable systems of subgroups
of $G$ (as the sets of minimal subgroups, cyclic subgroups or
proper terms of a composition series, respectively)?

\noindent{\bf Problem 3.5.} Given a finite group $G$, the
following function is well-defined
$$sd_G : L(G)\longrightarrow [0,1], \hspace{1mm} sd_G(H)=sd(H,G), \mbox{ for all } H\in
L(G).$$By Proposition 2.1, $sd_G$ is constant on each conjugacy
class of subgroups of $G$. Remark that the converse fails: take
the subgroups $H_1=<y>$ and $H_2=<xy>$ of $D_{2n}$; we have
$sd_G(H_1)=sd_G(H_2)=9/10$, but $H_1\nsim H_2$. Study other
properties of $sd_G$ (e.g. injectivity, monotony, ..., and so on),
as well as of the restriction of $sd_G$ to the set of conjugacy
classes of subgroups. Describe the finite groups $G$ for which
these functions satisfy certain conditions.

\noindent{\bf Problem 3.6.} Another interesting function can be
also associated to a finite group $G$, namely
$$sd(-,-) : L(G)\times L(G)\longrightarrow [0,1].$$Study this
function and its restrictions to some remarkable subsets of type
$L\times L$ of $L(G)\times L(G)$ (e.g. take $L=C(G)$, the poset of
cyclic subgroups of $G$). For an arbitrary $n\geq2$, generalize
the above function by defining
$$sd(\underbrace{-,-,...,-}_{n\ {\rm variables}}) : L(G)^n\longrightarrow [0,1],\hspace{2mm} sd(H_1,...,H_n)=$$
$$=\hspace{-0,5mm}\frax1{\dd\prod_{i=1}^n|L(H_i)|}\,|\{(K_1,...,K_n)\in\hspace{-1mm} \xmare{i=1}n\hspace{-0,5mm}L(H_i)\mid K_1\cdots K_n\hspace{-0,5mm}=\hspace{-0,5mm}K_{\sigma(1)}\cdots K_{\sigma(n)}, \forall\hspace{1mm} \sigma\in S_n\}|.$$

\noindent{\bf Acknowledgements.} The author is grateful to the
reviewer of \cite{4} for his suggestion to expand the list of open
problems on the subgroup commutativity degrees of finite groups in
a separate paper.

\end{document}